\documentclass{amsart}

\usepackage[all]{xy}
\usepackage{amssymb}
\usepackage{amsmath}
\usepackage{graphicx}
\usepackage[english]{babel}

\theoremstyle{definition}
\newtheorem{definition}{Definition}
\newtheorem{remark}[definition]{Remark}
\newtheorem{example}[definition]{Example}
\newtheorem{acknowledgement}{Acknowledgement}
\newtheorem{notation}[definition]{Notation}
\theoremstyle{plain}
\newtheorem{lemma}[definition]{Lemma}
\newtheorem{proposition}[definition]{Proposition}

\newtheorem{theorem}[definition]{Theorem}
\newtheorem{conjecture}[definition]{Conjecture}
\addtocounter{equation}{3}
\swapnumbers

\newcommand{\Vt}[1]{V^{\otimes #1}}
\def\Des{\operatorname{Des\!\!~}}
\def\Prim{\operatorname{Prim\ }}
\def\Ker{\operatorname{Ker\ }}
\def\Exp{\operatorname{E\!\!~ }}
\def\B{\operatorname{Ber\!\!~ }}
\def\exp{\operatorname{exp\!\!~ }}
\def\tr{\operatorname{tr\!\!~ }}

\def\log{\operatorname{log\!\!~ }}
\def\ad{\operatorname{ad\!\!~ }}

\def\Lie{\operatorname{Lie}}
\def\Id{\mathrm{Id}}
\def\Im{\mathrm{Im\ }}

\def\t{\otimes}
\def\KK{\mathbb K}
\def\QQ{\mathbb Q}
\def\NN{\mathbb N}
\def\RR{\mathbb R}
\def\CC{\mathbb C}
\def\HH{\mathcal H}
\def\gg{\mathfrak g}

\author[E.Burgunder]{Emily Burgunder}
\address{
Institut de Math\'ematiques et de mod\'elisation de Montpellier \\
UMR CNRS 5149\\
D\'epartement de math\'ematiques\\
Universit\'e Montpellier II\\
Place Eug\`ene Bataillon\\
34095 Montpellier CEDEX 5\\
France}
\email{burgunder@math.univ-montp2.fr}
\urladdr{www.math.univ-montp2.fr/{$\sim$}burgunder/}
\title{Eulerian idempotent and Kashiwara-Vergne conjecture}
\begin{document}
\begin{abstract}
By using  the interplay of the Eulerian idempotent and the Dynkin idempotent, we construct explicitly a particular symmetric solution $(F,G)$ of the first equation of the Kashiwara-Vergne conjecture:
$$x+y-\log(e^{y}e^{x})=(1-e^{-\ad x})F(x,y)+(e^{\ad y}-1)G(x,y) \ .$$
Then, we explicit all the solutions of the equation in the completion of the free Lie algebra generated by two indeterminates $x$ and $y$ thanks to the kernel of the Dynkin idempotent.
\end{abstract}
\maketitle
\let\languagename\relax
\section*{Introduction}
In \cite{KV}, M. Kashiwara and M.Vergne put forward a conjecture that implies the Duflo theorem on the local solvability of biinvariant differential operators on arbitrary finite Lie groups as well as a more general statement on convolution of invariant distributions:
\begin{conjecture}[Kashiwara-Vergne]\cite{KV}\label{KV}
For any Lie algebra $\gg$ of finite dimension, we can find Lie series $F$ and $G$ such that they satisfy:
\begin{enumerate}
\item $x+y-\log(e^{y}e^{x})=(1-e^{-\ad x})F(x,y)+(e^{\ad y}-1)G(x,y)$\label{KW1}.
\item $F$ and $G$ give $\gg$-valued convergent power series on $(x,y)\in \gg\times\gg$\label{KW2}.
\item \label{KW3}\begin{eqnarray*}
&&\tr(\ad x\circ\partial_x F;\gg)+\tr(\ad y\circ\partial_y G;\gg)=\\
&&\frac{1}{2}\tr(\frac{\ad x}{e^{\ad x} -1}+\frac{\ad y}{e^{\ad y} -1}+\frac{\ad \Phi(y,x)}{e^{\ad \Phi(y,x)} -1}-1;\gg)\ .
\end{eqnarray*}
\end{enumerate}
Here $\Phi(x,y)=\log(e^{x}e^{y})$ and $\partial_x F$ (resp. $\partial_y G$) is the $\mathbf{End}(\gg)$-valued real analytic function defined by
\begin{equation*}
\gg\ni a\mapsto \frac{d}{dt}F(x+ta,y)|_{t=0} \qquad \big( \textrm{resp.} \gg\ni a\mapsto \frac{d}{dt}G(x,y+ta)|_{t=0} \big)\ ,
\end{equation*}
and $\tr$ denotes the trace of an endomorphism of $\gg$.
\end{conjecture}
In the course of their proof they consider the first equation in the completion of the free Lie algebra generated by two indeterminates $x$ and $y$, denoted $\Lie(V)^\wedge$, where $V=\KK x\oplus \KK y$.

In this paper, we exhibit explicitly all the solutions of the first equation of the Kashiwara-Vergne conjecture in $\Lie(V)^\wedge$.

First, we display a particular solution thanks to a splitting of the equation and the use of two idempotents: the Eulerian idempotent $e$, and the Dynkin idempotent~$\gamma$.

The major problem to find solutions of equation (\ref{KW1}) is that there isn't any convenient basis  of the free Lie algebra that eases the calculation of the Baker-Campbell-Haussdorff series.
The Eulerian idempotent (cf.~\cite{L}) is the key to explicit the Baker-Campbell-Haussdorff series in terms of permutations.

We split the Baker-Campbell-Hausdorff series into two Lie formal power series $\Phi^-(y,x)$ and $\Phi^+(y,x)$, thanks to the Dynkin idempotent, such that they are in the image of $\exp(\ad(-x))-1$ (resp. image of $\exp(\ad y)-1$).

We split the equation (\ref{KW1}) into the sum of the two following equations:
\begin{eqnarray*}
\Phi^-(y,x)-x=(\exp\ad (-x)-1)F(x,y) \ , \\
\Phi^+(y,x)-y=(1-\exp\ad (y))G(x,y) \ ,
\end{eqnarray*}
and we give a unique  and explicit solution $(F_0(x,y),F_0(-y,-x))$ in terms of permutations. Hence we obtain a symmetric solution of equation (\ref{KW1}) by adding the two equations.

Any solution of equation (\ref{KW1}) on $\Lie(V)^\wedge$ is the sum of a particular solution and of a solution of the homogeneous equation:
\begin{equation*}
(\exp\ad (-x)-1)F(x,y)=(\exp\ad (y)-1)(G(x,y)) \ ,
\end{equation*}
for Lie formal power series $F$ and $G$.

We prove that any solution $(F,G)$ of the homogenous equation can be made explicit in terms of permutations and verifies that $xF(x,y)+yG(x,y)$ is in the kernel of the Dynkin idempotent, and conversely any element of the kernel of the Dynkin idempotent determines a solution of the homogeneous equation. 

Moreover these tools can be extended to prove a multilinearized version of the Kashiwara-Vergne conjecture.

There exists some solutions for some specific algebras by M. Kashiwara and M. Vergne \cite{KV}, F. Rouvi\`ere \cite{Ro} and M. Vergne \cite{V}. A. Alekseev and E. Meinrenken proved the existence of a solution in the general case in \cite{AM} using C. Torossian arguments  \cite{T} which are based on Kontsevich's work  cf.~\cite{K}. The latter has not  been made explicit and it is still unknown whether it is rational.

\medskip

The paper is organised as follows: in section \ref{sec:Notations} we set notations in section \ref{sec:Dynkin} and \ref{sec:Eulerian} we recall respectively the constructions of  the Dynkin idempotent and the Eulerian idempotent. Section \ref{sec:operator} is devoted to the study of  the operator $\exp\ad x -1$. The particular solution of the equation (\ref{KW1}) is constructed in section \ref{sec:particular}. This solution admits, in a certain way, a property of unicity which is treated in section \ref{sec:Unicity}. Then, we solve the homogeneous equation in section \ref{section:sanssecondmembre} and give a description of any solution of equation (\ref{KW1}) in the free Lie algebra. In section \ref{sec:otherview} we give another description of these solutions using another description of the kernel of the Dynkin idempotent due to F. Patras and C. Reutenauer. Section \ref{sec:multilinear} is devoted to the multilinear Kashiwara-Vergne conjecture.

\begin{acknowledgement}
I would like to thank Pr. J.-L. Loday and Pr. A. Brugui\`eres for their advisory, Dr. J.-M. Oudom for patient listening of some proofs, and Pr. M. Vergne for her careful reading of a preliminary version and for her useful remarks. 
\end{acknowledgement}

\bigskip

In this paper $\KK$ denotes a characteristic zero field, that is to say $\KK \supset\QQ$.
\section{Definitions and properties.}\label{sec:Notations}

We recall the definitions of bialgebras, convolution, tensor bialgebra and the free Lie algebra.

\subsection{Bialgebra and convolution}

A \emph{bialgebra} $\HH$ is a vector space endowed with an associative product $\mu:\HH\t\HH\rightarrow \HH$, a unit $u:\KK\rightarrow \HH$, a coassociative coproduct $\Delta:\HH\rightarrow \HH\t\HH$, and a counit $c:\HH\rightarrow \KK$ such that $\Delta$ is an algebra morphism or, equivalently, such that $\mu$ is a coalgebra morphism. 

The \emph{primitive part} of a bialgebra is the subvector space of $\HH$ defined as:
$$\Prim \HH :=\{x\in\HH:\Delta(x)=1\t x+x\t 1\} \ .$$

Let $f,g:\HH\longrightarrow \HH$ be two bialgebra morphisms. The \emph{convolution} of $f$ and $g$ is a bialgebra morphism defined as:
$$f\star g:=\mu\circ(f\t g)\circ\Delta|\HH\rightarrow\HH \ .$$
The convolution satisfies these easily verified propositions:
\begin{proposition}\label{starassociativity}
The convolution is associative.\hfill$\Box$
\end{proposition}
\begin{proposition}\label{starunity}
The convolution admits $u\circ c:\HH\rightarrow\HH $ for unit.\hfill$\Box$
\end{proposition}

\subsection{Tensor  bialgebra}\label{TV}

Let $V$ be a $\KK$-vector space. The \emph{tensor algebra} is the tensor module: 
$$T(V)= \KK\oplus V \oplus \Vt 2 \oplus\cdots\oplus\Vt n \oplus\cdots$$
endowed with the concatenation product $\mu:T(V)\t T(V)\rightarrow T(V):v_1\t\cdots\t v_n\bigotimes v_{n+1}\t\cdots\t v_{n+p}\mapsto v_1\t\cdots\t v_{n+p}$.
The tensor algebra can, moreover, be endowed with a unique coproduct (the co-shuffle)
$\Delta:T(V)\rightarrow T(V)\t T(V)$ such that $\Delta(v)= 1\t v+v\t 1$, making it into a cocommutative bialgebra.

\begin{remark} If the $\KK$-vector space $V$ is spanned by $< x_1,\ldots , x_n>$, then the tensor algebra $T(V)$ is spanned by all the tensors $x_{i_1}\t\cdots\t x_{i_m}$ where $i_j\in \{1,\ldots,n\}$, and $m\in\NN$. By the isomorphism induced by $x_{i_1}\t\cdots\t x_{i_m}\mapsto x_{i_1}x_{i_2}\cdots x_{i_m}$ for all $x_{i_j}\in  \{ x_1,\ldots, x_n\}$, the tensor algebra is isomorphic to the algebra of non-commutative polynomials in variables $x_1,\ldots,x_n$. 
\end{remark}
 
Following this remark we introduce the following notation in the algebra of non-commutative polynomials in variables $x_1,\ldots,x_n$:
\begin{definition}\label{def:x_part}
Let $p(x_1,\ldots,x_n)$ be a non-commutative polynomial in the variables $x_1,\ldots,x_n$. Then it determines uniquely $n$ polynomials $$b_1(x_1,\ldots,x_n),\ldots,b_n(x_1,\ldots,x_n)$$ such that $$p(x_1,\ldots,x_n)=x_1b_1(x_1,\ldots,x_n)+\cdots+x_n b_n(x_1,\ldots,x_n)\ .$$ 
We call $b_i(x_1,\ldots,x_n)$ the  \emph{$x_i$-part} of $p(x_1,\ldots,x_n)$ and denote it as $(p(x_1,\ldots,x_n))_{x_i}$.  
\end{definition}

Let $V$ be any vector space. The tensor bialgebra $T(V)$ verifies moreover this well-known connectedness property:
\begin{proposition}\label{connectedness}
Let $J:=\Id - u\circ c$. For any $w\in \Vt n$, $J^{\star n}(w)=0\ .$
\end{proposition}

\subsection{Action of the symmetric group}\label{subsec:symmgr}
Let $S_n$ be the symmetric group  acting on $\{1,\ldots,n\}$. It acts by the right on $\Vt n$ by permutation of the variables: for all $x_1,\ldots,x_n\in V$ and all $\sigma\in S_n$ the action is given by:
$$(x_1,\cdots,x_n)^\sigma=x_{\sigma(1)}\t\cdots\t x_{\sigma(n)}\ .$$

Let $\sigma\in S_n$ be a permutation. We denote by $\Des(\sigma)$ the set of \emph{descent} of $\sigma$  defined as  $\Des(\sigma):=\{i:\sigma(i)>\sigma(i+1)\}$ and by $d(\sigma)$
the \emph{number of descents} of $\sigma$ which is the number of integers $i$ such that $\sigma(i)> \sigma(i+1)$.

Denote by $D_{\{1,\ldots,k\}}$ the set of all the permutations $\sigma\in S_n$ such that its descent set $\Des(\sigma)$ is exactly the set $\{1,\ldots,k\}$.

\subsection{Free Lie algebra}

Let $V$ be a $\KK$-vector space. The \emph{free Lie algebra over $V$}, denoted $\Lie(V)$,
is defined by the following property:

any map $f : V\to \gg$, where $\gg$ is a Lie algebra,
extends uniquely to a Lie algebra morphism
$\tilde{f} : \Lie(V) \to \gg$. Diagrammatically, this would be read as the commutativity of the following diagram:
\begin{eqnarray*}
\forall\gg \ ,
&&\xymatrix{
V \ar@{^{(}->}[r]^{i} \ar[dr]_{\forall f}&\Lie(V) \ar[d]^{\exists !\tilde f}\\
&  \gg& \ .}
\end{eqnarray*}

It is known that $\Lie(V)$ can be identified with the subspace of $T(V)$ generated by $V$ under the bracket $[x, y]= x\t y-y\t x$. So we have: 
$$\Lie(V)= V\oplus[V,V]\oplus[V,[V,V]]\oplus\cdots\oplus\underbrace{[V,[V,[\cdots[V,V]\cdots]]]}_{n \textrm{ times}}\oplus\cdots \ .$$

\section{Free Lie algebra and Dynkin idempotent \cite{R,W}}\label{sec:Dynkin}

We recall the notion of Dynkin idempotent and some useful properties of the free Lie algebra.

\begin{definition} The \emph{Dynkin idempotent} is defined as 
$$\gamma:T(V)\rightarrow \Lie(V)\hookrightarrow T(V):v_1\t\cdots\t v_n\mapsto\frac{1}{n}[v_1,[v_2,[\cdots[v_{n-1},v_n]\cdots]]] \ .$$
The map $\gamma$ restricted to $\KK$ is null and its restriction to $V$ is the identity. It is clear from the definition of $\gamma$ and of $\Lie(V)$ that $\Im \gamma=\Lie(V)$.
We denote $\gamma|_{\Vt n}:\Vt n\rightarrow \Vt n$ by $\gamma_n:\Vt n\rightarrow \Vt n$ the restriction of the Dynkin idempotent on $\Vt n$.

We define the completion of the free Lie algebra $\Lie(V)^\wedge$ as $\Lie(V)^\wedge:=\Pi_{n\geq 0} \Lie(V)_n$.
\end{definition}

Moreover, this map satisfies the following properties.
\begin{proposition}\label{gamma}
For any $x,x_2,\cdots,x_n \in V$, with $n\geq 2$, one has: $$\gamma(x\t x_2\t \cdots\t x_n)=\frac{n-1}{n}ad(x)\gamma(x_2\t \cdots\t x_n)\ ,$$ where $ad(x):T(V)\rightarrow T(V):y\mapsto [x,y]$. 
\end{proposition}
\begin{proof}
Let $x,x_2,\cdots,x_n \in V$. A direct computation gives:
\begin{eqnarray*}
\gamma(x\t x_2\t \cdots\t x_n)=&\frac{1}{n}[x,[x_2,[ \cdots[x_{n-1}, x_n]\cdots]]]&\textrm{ by definition}\\
=&\frac{n-1}{n}[x,\gamma(x_2\t\cdots\t x_{n-1}\t x_n)] \ ,
\end{eqnarray*}
which completes the proof.
\end{proof}
\begin{proposition}[Friedrichs-Specht-Wever cf.~\cite{W}]
Suppose that $\KK$ is of characteristic zero. If $x\in\Vt n$, then the following are equivalent:
\begin{itemize}
\item $x\in \Lie(V)$,
\item $\Delta(x)=1\t x+x\t 1$,
\item $\gamma(x)=x$.
\end{itemize} 
\end{proposition}
Remark that, from this proposition, it  becomes clear that $\gamma$ is an idempotent that is to say $\gamma\circ\gamma=\gamma$.

Moreover this proposition gives another way to see $\Lie(V)$ as a subspace of $T(V)$ by:
\begin{proposition}\label{Prim} Let $V$ be a $\KK$-vector space. The primitive part of the tensor bialgebra is isomorphic to the free Lie algebra:
$\Prim T(V)\cong \Lie(V)\ .$
\end{proposition}

We state the following proposition with the notations of subsection \ref{subsec:symmgr}.
\begin{proposition}[Reutenauer cf.~\cite{R}]\label{explicitDynkin}
The restricted Dynkin idempotent $\gamma_n:\Vt n\rightarrow \Vt n$ can be made explicit in terms of permutations as follows: 
$$
\gamma_n(x_1\ldots x_n)=\frac{(-1)^{n-1}}{n}\sum_{k=0}^{n-1}(-1)^k\sum_{\sigma\in D_{\{1,\ldots,k\}}}(x_n\ldots x_1)^\sigma \ .
$$
\end{proposition}

\begin{proposition}\label{prop:KernelDynkin}
Let $V$ be a $\KK$-vector space.
For $n=0$ the kernel of the Dynkin idempotent restricted to $\KK$ is $\KK$. 

For $n=1$ the kernel of the Dynkin idempotent restricted to $V$ is $\{0\}$. 

Let $n\geq 2$. The kernel of the Dynkin idempotent restricted to $\Vt n$ is spanned by the elements $$(n-1)x_1\ldots x_n+\sum_{k=1}^{n-2}\sum_{\sigma\in D_{\{1,\cdots k\}}}(-1)^{n+k} (x_n\ldots x_1)^\sigma\ ,$$ for $x_i\in V$.
\end{proposition}

\begin{proof}
The cases $n=0$ and $n=1$ follow from the definition.
So, we focus on the case $n\geq 2$. As the Dynkin idempotent is a linear projector we can apply the usual trick which consists in writing any element $a\in T(V)$ as the following sum $(a-\gamma(a))+\gamma(a)$ belonging to $\Ker \gamma \oplus \Im \gamma$. 
Monomials $x_1\ldots x_n$ of degree $n$ are a basis for $\Vt n$, where $x_i\in V$. The kernel of the Dynkin idempotent are spanned by the elements $x_1\ldots x_n-\gamma(x_1\ldots x_n)$. 
To conclude it suffices to explicit the element $\gamma(x_1\ldots x_n)$ in terms of permutations thanks to the above proposition: indeed, the only permutation $\sigma$ of descent $\Des(\sigma)=\{1,\ldots,n-1\}$ is $\sigma=(n,n-1,\ldots,1)$ and  its sign is $(-1)^{2(n-1)}=1$.
\end{proof}

Adapting Patras and Reutenauer's description of the kernel of the Dynkin idempotent as a span (cf.~\cite{PR}) we get the following proposition:

\begin{proposition}\cite{C,PR}\label{kernelgamma}
The kernel of the Dynkin idempotent is spanned by 1 and the elements of the form $\gamma(a(x,y))a(x,y)$, where $a(x,y)$ is a non-commutative polynomial in the two indeterminates $x$ and $y$. 
\end{proposition}

\section{Eulerian idempotent \cite{L}}\label{sec:Eulerian}

We recall the notion of Eulerian idempotent and its link with the Baker-Campbell-Hausdorff series. Indeed, the Eulerian idempotent leads to an explicit formulation of the series in terms of permutations.

From now on, we consider $T(V)$ as a $\KK$-bialgebra with the concatenation $\mu$, a unit $u$, the co-shuffle $\Delta$, and a counit $c$. The convolution product is denoted by $\star$. We will say Lie series for formal power Lie series.

\begin{definition} Define the map $J:=Id-u\circ c:T(V)\rightarrow T(V)$. The \emph{Eulerian idempotent} $e$ is 
the following endomorphism of $T(V)$ defined by the formal power series: 
$$e:=\log^\star(Id)=\log^\star(J+u\circ c)=J-\frac{J^{\star 2}}{2}+\frac{J^{\star 3}}{3}+\cdots+(-1)^{k-1}\frac{J^{\star k}}{k}+\cdots\ ,$$
where $J^{\star n}=\underbrace{J\star J\star\cdots\star J}_{n \textrm{ times}}$ .
\end{definition}
 Proposition \ref{starassociativity} assures that $J^{\star n}$ is well-defined. The map $J:T(V)\rightarrow T(V)$ is the identity on $V\oplus V^{\t 2}\oplus\cdots\oplus V^{\t n}\oplus\cdots$ and is null on $\KK$. Thanks to the connectedness property of proposition \ref{connectedness}, the restriction of $e$ to $\Vt n$ is polynomial and is equal to:
 $$e|_{\Vt n} =J-\frac{J^{\star 2}}{2}+\frac{J^{\star 3}}{3}+\ldots+(-1)^n\frac{J^{\star n-1}}{n-1}\ . $$
The restriction of the map $e$ to $\Vt n$ is denoted $e_n:\Vt n\rightarrow\Vt n$.
\begin{proposition}\cite{L}
The map $e:T(V)\rightarrow T(V)$ verifies the following properties:
\begin{itemize}
\item $\Im e=\Lie(V)$,
\item the map $e$  is an idempotent, i.e. $e\circ e=e$.
\end{itemize}
\end{proposition}

The Baker-Campbell-Hausdorff series is a formal power series $\Phi(x,y)=\sum_{n\geq 1} \Phi_n(x,y)$, where $\Phi_n(x,y)$ is a homogeneous polynomial of degree $n$, in non-commutative variables $x$ and $y$, defined by the equation:
\begin{equation}
\exp(x)\exp(y)=\exp(\Phi(x,y))\ ,\label{BCHeq}
\end{equation}
where $\exp$ denotes the exponential series.

This formal power Lie series can be extended to $n$ variables, by defining $\Phi(x_1,\cdots,x_n)$ as: 
$$\exp(x_1)\cdots \exp(x_n)=\exp(\Phi(x_1,\ldots,x_n)) \ .$$
Let $\Phi_m(x_1,\ldots,x_n)$ denote the homogeneous part of $\Phi(x_1,\ldots,x_n)$ of total degree $m$, and $\varphi_n(x_1,\ldots,x_n)$ the multilinear part of $\Phi_n(x_1,\ldots,x_n)$, (replace $x_i^2$ by $0$, for all $i$).
\begin{proposition}[Dynkin cf.~\cite{L}]\label{Dynkin}
The following equality holds:
$$\Phi_m(x_1,\cdots,x_n)=\sum_{\begin{array}{c}
\scriptstyle i_1+\cdots+i_n=m \\
\scriptstyle i_j\geq 0
\end{array}}\frac{1}{i_1!\ldots i_n!}\varphi_m(\underbrace{x_1\ldots,x_1}_{i_1},\dots,\underbrace{x_n\ldots,x_n}_{i_n})\ .$$
\end{proposition}

The next proposition relates the Baker-Campbell-Hausdorff series to the Eulerian idempotent.
\begin{proposition}\cite{L}\label{relates}
The following equality holds: 
$\varphi_n(x_1,\ldots,x_n)=e_n(x_1,\cdots,x_n)$.
\end{proposition} 

\subsection{Explicitation of $e_n$}\label{Explicit} 

The Eulerian idempotent $e_n:\Vt n\rightarrow \Vt n$ can be made explicit as a formal power series of permutations with the notations of subsection \ref{subsec:symmgr}:
\begin{proposition}\cite{L}\label{explicit}The Eulerian idempotent has the following explicit formula in terms of permutations: $e_n=\sum_{\sigma \in \mathcal S_n}c_\sigma (.)^\sigma$, where 
$c_\sigma=(-1)^{d(\sigma)}\binom{n-1}{d(\sigma)}^{-1}$.
\end{proposition} 
\noindent Here $\binom{n}{p}$ denotes the binomial number.

If we restrict ourselves to $T(V)$, with $V=\KK x\oplus \KK y$,
then the above formulas lead to: 
\begin{proposition}\label{Vspecialise} Let $V=\KK x\oplus \KK y$. The Baker-Campbell-Hausdorff series can be made explicit in terms of permutations:
\begin{equation*}\Phi(x,y)=\sum_{n\geq 1}\sum_{\begin{array}{c}
\scriptstyle i+j=n \\
\scriptstyle i,j\geq 1
\end{array}}\frac{1}{i!}\frac{1}{j!}\sum_{\sigma \in \mathcal S_n}c_\sigma (\underbrace{x,\ldots,x}_{i},\underbrace{y,\ldots,y}_{j})^\sigma \ ,
\end{equation*}
where 
$c_\sigma=(-1)^{d(\sigma)}\binom{n-1}{d(\sigma)}^{-1}$.
\end{proposition}
\begin{proof}
By definition $\Phi(x,y)=\sum_{n\geq 1}\Phi_n(x,y)$. Restricting propositions \ref{relates} and \ref{Dynkin}  to the two variables $x$ and $y$, we get that:
\begin{equation*}\Phi(x,y)=\sum_{n\geq 1}\sum_{\begin{array}{c}
\scriptstyle i+j=n \\
\scriptstyle i,j\geq 1
\end{array}}\frac{1}{i!}\frac{1}{j!}e_n(\underbrace{x,\ldots,x}_{i},\underbrace{y,\ldots,y}_{j}) \ .
\end{equation*}
Proposition \ref{explicit} completes the proof.  
\end{proof}

\subsection{Symmetric properties of the Eulerian idempotent}
Let $\omega$ denote the following permutation $\omega(1,\ldots,n):=(n,n-1,\ldots,2,1)$.
With the above explicitation of the Eulerian idempotent, we remark the following symmetry:
\begin{proposition}\label{ensymmetry}
The Eulerian idempotent verifies the following symmetry: $e_n=(-1)^{n+1}e_n\circ()^\omega$.
\end{proposition} 
\begin{proof}
We first note that $(-1)^{n+1}a_{\sigma}=a_{ \omega\circ \sigma}$. Indeed, we have:
\begin{eqnarray*}
a_{\omega\circ \sigma }&=&(-1)^{d(\omega\circ \sigma )}\binom{n-1}{d(\omega\circ \sigma )}^{-1} 
=(-1)^{n-1}(-1)^{d(\sigma)}\binom{n-1}{n-1-d(\sigma)}^{-1}
\\
&=&(-1)^{n-1}(-1)^{d(\sigma)}\binom{n-1}{d(\sigma)}^{-1} 
=(-1)^{n-1}a_{\sigma},
\end{eqnarray*}
as $d(\omega\circ \sigma)=n-1-d(\sigma)$.
This gives the expected property.

As a consequence we have:
$$ e_n=\sum_{\sigma \in \mathcal S_n}c_\sigma ()^\sigma=\sum_{\omega\circ \sigma \in \mathcal S_n} 
a_{\omega\circ \sigma } ()^{\omega\circ \sigma}  = (-1)^{n-1}\sum_{\sigma\in \mathcal S_n} 
a_{\sigma} ()^{\omega\circ \sigma} =(-1)^{n+1}e_n\circ()^\omega  .$$
\end{proof}

As a consequence we can prove that: 
\begin{lemma}\label{lemma:Asymmetry}
The $x$-part of the Eulerian idempotent verifies the following  symmetry property:
$$(e_{n}(\underbrace{x,\ldots,x}_{i},\underbrace{y,\ldots,y}_{j})_x)=(e_{n}(\underbrace{-y,\ldots,-y}_{j},\underbrace{-x,\ldots,-x}_{i})_y)$$
\end{lemma}
\begin{proof}Let $\omega$ denote the following permutation $\omega(1,\ldots,n):=(n,n-1,\ldots,2,1)$.
Proposition \ref{ensymmetry} gives a symmetry property of the Eulerian idempotent: $e_n=(-1)^{n+1}e_n\circ()^\omega$.
From this symmetry property  we deduce the following:
\begin{eqnarray*}
&&x(e_{n-1}(\underbrace{x,\ldots,x}_{i-1},\underbrace{y,\ldots,y}_{j})_x) +y (e_{n-1}(\underbrace{x,\ldots,x}_{i},\underbrace{y,\ldots,y}_{j-1})_y)\\ &&\qquad=(-1)^{n-1}y(e_{n-1}(\underbrace{y,\ldots,y}_{j-1},\underbrace{x,\ldots,x}_{i})_x)+ (-1)^{n-1}x(e_{n-1}(\underbrace{y,\ldots,y}_{j},\underbrace{x,\ldots,x}_{i-1})_y)\\
&&\qquad=y (e_{n-1}(\underbrace{-y,\ldots,-y}_{j-1},\underbrace{-x,\ldots,-x}_{i})_x)+x(e_{n-1}(\underbrace{-y,\ldots,-y}_{j},\underbrace{-x,\ldots,-x}_{i-1})_y)\ .
\end{eqnarray*}

Therefore, we get $(e_{n}(\underbrace{x,\ldots,x}_{i},\underbrace{y,\ldots,y}_{j})_x)=(e_{n}(\underbrace{-y,\ldots,-y}_{j},\underbrace{-x,\ldots,-x}_{i})_y)$, which ends the proof.
\end{proof}

\section{The operator $\Exp( x)$}\label{sec:operator}
This section is devoted to a study of the image and the kernel of the operator $\exp\ad x-1:T(V)\rightarrow T(V)$, restricted to $\exp\ad x-1:\Lie(V)\rightarrow \Lie(V)$, which is needed in severals proofs.

Let $\KK$ be a characteristic zero field, i.e. $\KK \supset\QQ$. From now on, $V$ is the following $\KK$-vector space $V=\KK x\oplus \KK y$, and $T(V)$ the tensor bialgebra defined in sections \ref{TV} and \ref{sec:Eulerian}. Recall that this tensor bialgebra is isomorphic to the non-commutative polynomial bialgebra in two variables $x$ and $y$. 
For conveniency, we will say Lie series for Lie formal power series.

We will use the following notation:
\begin{notation}\label{notation}
Denote $\Exp( x):\Lie(V)\rightarrow\Lie(V)$ the map $\Exp(x):= \exp (\ad x) -1=\sum_{n\geq 1}\frac{(\ad x)^n}{n!}$.
\end{notation}

\begin{notation}\label{not:Bernoulli}
Denote $B_n$ be the $n$-th Bernoulli number defined as:
$$B_n:=(-1)^n\sum_{k=1}^{2n+1}\frac{(-1)^k}{k}\binom{2n+1}{k}\sum_{r=0}^{k}r^{
2n}\ .$$

\end{notation}

\begin{notation}\label{not:SerieBernoulli}
Denote $\B( x):=\frac{ad(x)}{\Exp((x))}:\Lie(V)\rightarrow \Lie(V)$ the map defined as: $$\B( x):=\sum_{n\geq 0}B_k \frac{\ad(x)^k}{k!}\ .$$
\end{notation}

\begin{proposition}\label{prop:inverse}
The maps $\B(x),\Exp(x):T(V)\rightarrow T(V)$ verify the property: 
$$\B(x)\circ\Exp(x)=\Exp(x)\circ\B(x)=\ad (x) \ .$$
\end{proposition}

\begin{proof}
The composition of series of operators is the multiplication of the series.
As series $\frac{t}{\exp(t)-1}$ defined as $\sum_{n\geq 0}B_k \frac{t^k}{k!}$ admits the property: 
$$(\exp(t)-1) \frac{t}{\exp(t)-1}= \frac{t}{\exp(t)-1}(\exp(t)-1)=t\ ,$$ this ends the proof. 
\end{proof}

\begin{proposition}\label{Image}
The image of $\Exp (  x) : \Lie(V)\to \Lie(V)$ is:
$$\Im \Exp(  x) = \Im \ad x \ .$$
\end{proposition}

\begin{proof}
Clearly $\Im \Exp(  x) \subset \Im \ad x$. Moreover, any element in $\Im \ad x$ is of the form $[x,\alpha]\in\Im\ad x$, where $\alpha\in \Lie(V)$. Then, $\beta:=\B( x)(\alpha)=\sum_{n\geq 0} \frac{B_{n}}{n!} (\ad x)^n (\alpha)\in \Lie(V)$ is such that $$\Exp(  x) (\beta)=\Exp ( x)\circ(\B  x)(\beta)=\ad x(\alpha)=[x,\alpha]\ .$$ 
\end{proof}

\begin{lemma}\label{Kernelad}
The kernel of the operator $\ad x:T(V)\rightarrow T(V)$ is spanned by the elements $x^n$:
$$
\Ker \ad x=\KK[x]\ .
$$
\end{lemma}

\begin{proof}
Let $p(x,y)$ be a non-commutative polynomial in the kernel of $\ad x$, i.e. such that $\ad x\ p(x,y)=0$. This non-commutative polynomial can be decomposed as the sum of its homogeneous part: $p(x,y)=\sum_{n\geq 0}p_n(x,y)$, and each $p_n(x,y)$ verifies $\ad x \ p_n(x,y)=0$.

The proof is based on induction on the total degree $n$ of the homogeneous polynomial.

Let $n=0$, then $p_0(x,y)=\lambda$, where $\lambda\in\KK$, verifies the equation. 

Let $n=1$, then a generic non-commutative  homogeneous polynomial of degree one can be made explicit on the basis of $V$ as $p_1(x,y)=\lambda x+\mu y$, where $\lambda,\mu\in\KK$. Then, in order to verify the equation $\ad x\ p_n(x,y)=0$, we must take  $\mu=0$ and $\lambda\in\KK$.

For the degree $n$, suppose by induction that  $p_k(x,y)=\lambda x^k$, for any $k\leq n-1$.  The homogeneous non-commutative polynomial $p_n(x,y)$ can be split into:  
$$p_n(x,y)=xq_{n-1}(x,y)+yr_{n-1}(x,y)\ ,$$
 where $q_{n-1}(x,y)$, $r_{n-1}(x,y)$ are respectively the  $x$-part and the  $y$-part of $p_n(x,y)$ , i.e. non-commutative polynomials in variables $x$ and $y$ of total degree $n-1$ verifying the above equation. And so the equation $\ad x\ p_n(x,y)=0$ can be rewritten as:
\begin{eqnarray*}
&xxq_{n-1}(x,y)-xq_{n-1}(x,y)x+xyr_{n-1}(x,y)-yr_{n-1}(x,y)x&=0\\
\textrm{or, }&x(xq_{n-1}(x,y)-q_{n-1}(x,y)x+yr_{n-1}(x,y))-yr_{n-1}(x,y)x&=0 \ .
\end{eqnarray*}
As the above identity is the nullity of a non-commutative polynomial, by identification we have that $r_{n-1}(x,y)=0$ and we are left with the equation:
\begin{eqnarray*}
&xq_{n-1}(x,y)-q_{n-1}(x,y)x&=0\\
\textrm{or, }&\ad(x)q_{n-1}(x,y)&=0\ 
\end{eqnarray*}
to solve. Applying the induction hypothesis gives that $p_n(x,y)=\lambda x^n$, where $\lambda\in\KK$, and it ends the proof.
\end{proof}

\begin{lemma}\label{KerneladLie}
The kernel of the operator $\ad x:\Lie(V)\rightarrow \Lie(V)$ is:
$$
\Ker \ad x=\KK x\ .
$$
\end{lemma}

\begin{proof}
As the kernel of $\ad x:T(V)\rightarrow T(V)$ is the polynomial algebra in indeterminate $x$, the kernel of $\ad x:\Lie(V)\rightarrow \Lie(V)$ is $\Lie(V)\cap\KK[x]=\KK x$. This concludes the proof.
\end{proof}

\begin{lemma}\label{lemma:passeadexp}
Let $S(x,y)$ and $P(x,y)$ be two Lie polynomials. If these polynomials verify:
\begin{equation}\label{eq:passage}
\Exp ( -x)S(x,y)=\ad (-x)P(x,y) \ ,
\end{equation}
then,
$$S(x,y)=\B(-x)P(x,y)+\lambda x\ ,$$
for a certain $\lambda\in\KK$.
\end{lemma}
\begin{proof}
The fact that the pair $(S(x,y)=\B(-x)P(x,y)+\lambda x, P(x,y))$ verifies equation (\ref{eq:passage}) is straigthforward.
Conversely, suppose that the pair $(S(x,y),P(x,y))$ verifies equation (\ref{eq:passage}):
$$
\Exp (-x)S(x,y)=\ad (-x)P(x,y) \ .
$$
Then, multiplying on the left  side by the operator $\B(-x)$ leads to:
$$
\ad (-x)S(x,y)=\B(-x)\ad(-x)P(x,y)\ .
$$
The following identity follows from the commutativity of $\ad(-x)$ and $\Exp(-x)$~:
$$
\ad (-x)S(x,y)=\ad (-x)\B(-x)P(x,y)\ ,
$$
therefore the Lie  series $S(x,y)-\B(-x)P(x,y)$ is in the kernel of the morphism $\ad(-x)$ exhibited in proposition \ref{KerneladLie}. This ends the proof as there exists $\lambda\in\KK$ such that:
$$S(x,y)=\B(-x)P(x,y)+\lambda x\ .$$
\end{proof}

\begin{proposition}\label{Kernel}
The kernel of the operator $\Exp (  x) : \Lie(V)\to \Lie(V)$ is:
$$\Ker \Exp(  x) =\KK x\ .$$
\end{proposition}

\begin{proof}
By proposition \ref{Image} for any Lie polynomial $P(x,y)$ there exists a Lie polynomial $p(x,y)$ such that 
\begin{equation}\label{eq:kerexp}
\Exp(x)P(x,y)=ad(x)p(x,y)\ .
\end{equation}
By proposition \ref{lemma:passeadexp} there exists $\mu\in\KK$ such that:
$$
P(x,y)=\B( x)p(x,y) +\mu x\ .
$$
Remark that to determine the kernel of the operator $\Exp(x)$ we will use the kernel of the adjunction, i.e. if $P(x,y)\in\Ker\Exp((x))$ then we have $p(x,y)\in\Ker\ad(x)$ and by proposition \ref{Kernelad} we conclude that $p(x,y)=\lambda x$
and therefore we have:
$$
P(x,y)=\B( x)(\lambda x) +\mu x\ .
$$
By expanding the inverse operator we have,
\begin{eqnarray*}
P(x,y)&=&\sum_{k\geq 0}\frac{B_k}{k!} \ad(k)^{n}( \lambda x) + \mu x\\
&=&\frac{B_0}{0!}(\lambda x) + \mu x\\
&=&(\lambda+\mu)x\ ,
\end{eqnarray*}
which concludes the proof.
\end{proof}

\section{A particular solution of equation \ref{KW1}}\label{sec:particular}

In this section, we construct an explicit symmetric solution to the equation (\ref{KW1}) of the Kashiwara-Vergne conjecture by taking the Dynkin idempotent of a split of the Eulerian idempotent.

We adopt notation  \ref{notation}
and use the Baker-Campbell-Hausdorff series (proposition \ref{Vspecialise}) to rewrite the Kashiwara-Vergne first equation.
\begin{proposition}\label{rewrite}
The first equation of the Kashiwara-Vergne conjecture (\ref{KW1}) is equivalent to:
\begin{equation}
\sum_{n\geq 2}\Phi_n(y,x)=\Exp(-x)F(x,y)-\Exp( y) G(x,y) \ .\label{KW1modif}
\end{equation}
\end{proposition}
\begin{proof}
By equation (\ref{BCHeq}), $\log(\exp y \exp x)=\log\exp\Phi(y,x)=\Phi(y,x)=\sum_{n\geq 1}\Phi_n(y,x)$. Moreover, by proposition \ref{explicit}, $\Phi_1(x,y)=e_1(x)+e_1(y)=x+y$. So the left part of the Kashiwara-Vergne first equation,
$$x+y-\log(e^{y}e^{x})=(1-e^{-\ad x})F(x,y)+(e^{\ad y}-1)G(x,y)\ ,$$ can be rewritten as $\sum_{n\geq 2}\Phi_n(y,x)=(\exp (\ad -x)-1)F(x,y)-(\exp(\ad y)-1)G(x,y)$. The use of the notation \ref{notation} completes the proof.  
\end{proof}

With the notations of section \ref{Explicit}, we define the two following polynomials:
\begin{definition}\label{phihalf}
Define $\Phi_x(x,y)$ as the sum of all the monomials of $\sum_{n\geq 2}\Phi_n(x,y)$ beginning with the indeterminate $x$:
$$\Phi_x(x,y):=x(\sum_{n\geq 2}(\Phi_n(x,y))_x)\ ,$$
with the notation of definition \ref{def:x_part}.

We denote $\Phi^+(x,y):=\gamma(\Phi_x(x,y))$ and  $\Phi^-(x,y):=\gamma(\Phi_y(x,y))$.
\end{definition}
Remark that $\Phi_x(x,y)$ and $\Phi_y(x,y)$ are non-commutative  series and not Lie  series. Taking their Dynkin idempotent forces them to be Lie  series and to verify the following proposition :

\begin{proposition}\label{prop:welldef}
The formal power series $\Phi^+(x,y)$ and $\Phi^-(x,y)$ satisfy the following property: 
$$\Phi^+(y,x)\in\Im\Exp(y)\ ,\qquad \Phi^-(y,x)\in\Im\Exp(-x)\ .$$
\end{proposition}
\begin{proof}
By definition \ref{phihalf} the Lie series  $\Phi^+(x,y)$ is defined as $\Phi^+(y,x)=\gamma(y(\Phi(y,x))_y)$. Proposition \ref{gamma} assures that this Lie series is in the image of $\ad(y)$. Then applying proposition \ref{Image} ends the proof. The other property is proved analogously. 
\end{proof}

These formal power series can be made explicit in terms of permutations as follows:
\begin{proposition} 
The formal power series defined above have the following expression in terms of permutations:
\begin{eqnarray*}
\Phi^+(x,y):=\sum_{n\geq 2}\sum_{\begin{array}{c}
\scriptstyle i+j=n \\
\scriptstyle i,j\geq 1
\end{array}} \frac{1}{i!}\frac{1}{j!}\sum_{\begin{array}{c}
\scriptstyle \sigma \in S_n \\
\scriptstyle \sigma^{-1}(1)\in\{1,\ldots, i\}
\end{array}}\gamma\circ c_\sigma (\underbrace{x,\ldots,x}_{i},\underbrace{y,\ldots,y}_{j})^\sigma \ ,\\
\Phi^-(x,y):=\sum_{n\geq 2}\sum_{\begin{array}{c}
\scriptstyle i+j=n \\
\scriptstyle i,j\geq 1
\end{array}} \frac{1}{i!}\frac{1}{j!}\sum_{\begin{array}{c}
\scriptstyle \sigma \in S_n \\
\scriptstyle \sigma^{-1}(1)\in\{i+1,\ldots, n\}
\end{array}}\gamma\circ c_\sigma (\underbrace{x,\ldots,x}_{i},\underbrace{y,\ldots,y}_{j})^\sigma \ ,
\end{eqnarray*}
where $c_\sigma=(-1)^{d(\sigma)}\binom{n-1}{d(\sigma)}^{-1}$.
\end{proposition}
\begin{proof}
Remark that $\sum_{n\geq 2}\Phi(x,y)$ can be made explicit thanks to its link with the Eulerian idempotent (cf. proposition \ref{Vspecialise}). Taking its $x$-part  is to restrict the explicit version of $\sum_{n\geq 2}\Phi(x,y)$ to the permutations $\sigma$ such that $\sigma^{-1}(1)\in\{1,\ldots,i\}$ which garantees that the monomial will start with an $x$.
\end{proof}
We could make this formula even more explicit by using proposition \ref{explicitDynkin}.

These polynomials split the left part of the Kashiwara-Vergne conjecture (\ref{KW1modif}):
\begin{proposition}\label{cut}
The two formal power series defined in definition \ref{phihalf} verify the property:
$$\Phi_n(x,y)=\Phi_n^+(x,y)+\Phi_n^-(x,y)\ .$$
\end{proposition}
\begin{proof}
By definition of $\Phi_n^+(x,y)$ and $\Phi_n^-(x,y)$, we have:
\begin{eqnarray*}
\Phi_n^+(x,y)+\Phi_n^-(x,y)&=&
\gamma(x((\Phi_n(x,y))_x))+\gamma(y(\Phi_n(x,y)_y))\ ,
\end{eqnarray*}
which is equal to:
$$\Phi_n^+(x,y)+\Phi_n^-(x,y)= \gamma(\Phi_n(x,y))=\Phi_n(x,y)\ ,$$ 
as $\Phi_n(x,y)$ is a Lie polynomial.
\end{proof}

This splitting admits moreover a certain symmetry:
\begin{lemma}\label{phisymmetry}
The split Baker-Campbell-Hausdorff series verify the following anti-symmetric property: $\Phi^+(x,y)=-\Phi^-(-y,-x)\ .$
\end{lemma}
\begin{proof}
We use the definition $\Phi_x(x,y)=x(\sum_{n\geq 2}(\Phi_n(x,y))_x)$, to prove the symmetry property. By proposition \ref{Vspecialise} the formal power series $\sum_{n\geq 2}\Phi(x,y)$ can be made explicit thanks to its link with the Eulerian idempotent: \begin{equation*}\Phi^+(x,y)=\sum_{n\geq 2}\sum_{\begin{array}{c}
\scriptstyle i+j=n \\
\scriptstyle i,j\geq 1
\end{array}}\frac{1}{i!}\frac{1}{j!}\gamma(x(e_n(\underbrace{x,\ldots,x}_{i},\underbrace{y,\ldots,y}_{j}))_x) \ .
\end{equation*}
By proposition \ref{lemma:Asymmetry}, we have:
\begin{equation*}\Phi^+(x,y)=\sum_{n\geq 2}\sum_{\begin{array}{c}
\scriptstyle i+j=n \\
\scriptstyle i,j\geq 1
\end{array}}\frac{1}{i!}\frac{1}{j!}\gamma(x(e_n(\underbrace{-y,\ldots,-y}_{j},\underbrace{-x,\ldots,-x}_{i}))_y) \ .
\end{equation*}
And the symmetry is proven as:
 \begin{equation*}-\Phi^-(-y,-x)=\sum_{n\geq 2}\sum_{\begin{array}{c}
\scriptstyle i+j=n \\
\scriptstyle i,j\geq 1
\end{array}}\frac{1}{i!}\frac{1}{j!}\gamma(x(e_n(\underbrace{-y,\ldots,-y}_{i},\underbrace{-x,\ldots,-x}_{j}))_y) \ .
\end{equation*} 
\end{proof}

In order to simplify the particular solution of equation (\ref{KW1}) we construct the following polynomial:
\begin{definition}\label{def:AetB}
Define the   Lie series $a(x,y)$ as follows:
\begin{eqnarray*}
a(x,y)&:=&\sum_{n\geq 1}\frac{n}{n+1}\sum_{\begin{array}{c}
\scriptstyle i+j=n \\
\scriptstyle i,j\geq 1
\end{array}} \frac{1}{(i+1)!}\frac{1}{j!}\gamma( (e_n(\underbrace{x,\ldots,x}_{i},\underbrace{y,\ldots,y}_{j}))_x) \ .
\end{eqnarray*}
\end{definition}
As before, though $e_n$ is a Lie polynomial, its $x$-part is not a Lie polynomials in general. So, in order to force $a(x,y)$ to become a Lie  series, we have applied the Dynkin idempotent $\gamma$ to $e_n(x,\ldots,x,y,\ldots,y)_x$.

Proposition \ref{explicit} gives an explicit version of the Eulerian idempotent $e:T(V)\rightarrow T(V)$ which permits us to define explicitly these two Lie polynomials: 
\begin{proposition}\label{prop:en_xenperm} 
Let $\sigma\in S_n$, we denote $\tilde\sigma$ the image of the $n-1$ last variables: $\sigma(1\ldots n)=(\sigma(1),\tilde\sigma(2,\ldots,n))$.
Then, $(e_n(\underbrace{x,\ldots,x}_{i},\underbrace{y,\ldots,y}_{j}))_x$ can be made explicit in the following way:
\begin{eqnarray*}
&&(e_n(\underbrace{x,\ldots,x}_{i},\underbrace{y,\ldots,y}_{j}))_x=\sum_{\begin{array}{c}
\scriptstyle \sigma\in S_n \\
\scriptstyle \sigma^{-1}(1)\in\{1,\ldots, i\}
\end{array}}c_\sigma \gamma\circ(\underbrace{x,\ldots,x}_{i-1},\underbrace{y,\ldots,y}_{j})^{\tilde\sigma}\ .
\end{eqnarray*}\hfill$\Box$
\end{proposition}
This formula could be made more explicit by proposition \ref{explicitDynkin}.

Now we can give the definition of the particular solution of the equation (\ref{KW1}):
\begin{definition}\label{def:FetG}
Let $\KK$ be a  characteristic zero field, and $V$ the vector space defined by $V=\KK x\oplus \KK y$. 
Let $a(x,y)$ be the  Lie  series defined in definition \ref{def:AetB}.
We
define the Lie  series $F_0(x,y)$  as:
\begin{eqnarray*}
&&F_0(x,y):=-\sum_{n\geq 0}\frac{B_n}{n!}(-1)^n(\ad x)^{ n}\circ a(-x,-y)\ .
\end{eqnarray*}
\end{definition}
Note that the  Lie series is well defined as restricted to elements of degree $n$, $a(-x,-y)$  is polynomial. Therefore it is polynomial when restricted to a degree $n$.

Thanks to the two polynomials $\Phi^+(x,y)$ and $\Phi^-(x,y)$ defined from the Baker-Campbell-Hausdorff series (cf. definition \ref{phihalf}) we split the Kashiwara-Vergne first equation into the following equation:
\begin{proposition}\label{splitKV}
The equation:
\begin{equation}
\Phi^-(y,x)=\Exp(-x)F(x,y) \ . \label{splitequation}
\end{equation}
admits the Lie  series $F_0(x,y)$ defined in definition \ref{def:FetG} as solution on $\Lie(V)$.
\end{proposition}

\begin{proof}\label{demo}
The equation \ref{splitequation} is well-defined by proposition \ref{prop:welldef}. 

By notation \ref{not:Bernoulli} the solution can be rewritten as $F_0(x,y)=-\B (-x)\circ a(-x,-y)$. It verifies the following equalities by proposition \ref{prop:inverse}, definition \ref{phihalf}, proposition \ref{gamma} and proposition \ref{lemma:Asymmetry} respectively :
\begin{eqnarray*}
\Exp(-x)F_0(x,y)&=&-\Exp(-x)\B (-x)\circ a(-x,-y)\\
&=&\ad x\circ a(-x,-y)\\
&=&\ad x \circ \sum_{n\geq 2}\frac{n-1}{n}\sum_{\begin{array}{c}
\scriptstyle i+j=n \\
\scriptstyle i\geq 2,j\geq 1
\end{array}} \frac{1}{i!}\frac{1}{j!}\gamma\circ (-1)^{n-1}(e_{n-1}(\underbrace{x,\ldots,x}_{i-1},\underbrace{y,\ldots,y}_{j})_x)\\
&=&\sum_{n\geq 2}\sum_{\begin{array}{c}
\scriptstyle i+j=n \\
\scriptstyle i \geq 2,j\geq 1
\end{array}} \frac{1}{i!}\frac{1}{j!}(-1)^{n-1}\gamma( x (e_{n-1}(\underbrace{x,\ldots,x}_{i-1},\underbrace{y,\ldots,y}_{j})_x))\\
&=&\sum_{n\geq 2}\sum_{\begin{array}{c}
\scriptstyle i+j=n \\
\scriptstyle i \geq 2,j\geq 1
\end{array}} \frac{1}{i!}\frac{1}{j!}\gamma( x (e_{n-1}(\underbrace{y,\ldots,y}_{j},\underbrace{x,\ldots,x}_{i-1})_y))\\
&=& \sum_{n\geq 2}\sum_{\begin{array}{c}
\scriptstyle i+j=n \\
\scriptstyle i,j\geq 1
\end{array}} \frac{1}{i!}\frac{1}{j!}\sum_{\begin{array}{c}
\scriptstyle \sigma \in S_{n} \\
\scriptstyle \sigma^{-1}\in\{i+1,\cdots, n\}
\end{array}}\gamma\circ c_\sigma (\underbrace{y,\cdots,y}_{i},\underbrace{x,\cdots,x}_{j})^\sigma\\
&=& \Phi^-(y,x)\ .
\end{eqnarray*} 
This ends the proof.
\end{proof}

Thanks to  lemma \ref{phisymmetry} we can prove:
\begin{proposition}\label{prop:equationsplitsymmetry}
Let $F(x,y)$ be a Lie series which is a  solution of the split equation:
$$\Phi^-(y,x)=\Exp(- x)F(x,y) \ . $$
Then, $F(-y,-x)$ is a solution of:
$$\Phi^+(y,x)=-\Exp ( y)G(x,y)\ .$$
\end{proposition}

\begin{proof}
Let $F(x,y)$ be a Lie polynomial solution of the split equation $\Phi^-(y,x)=\Exp(- x)F(x,y)$. By exchanging $x$ and $-y$ we have $\Phi^-(-x,-y)=\Exp( y)F(-y,-x)$. Then by lemma  \ref{phisymmetry} we get $\Phi^-(-x,-y)=\Phi^+(y,x)$, and so $\Phi^+(y,x)=\Exp( y)F(-y,-x)$. Then $G(x,y)=F(-y,-x)$ is a solution of $\Phi^+(y,x)=-\Exp ( y)G(x,y)$, which ends the proof.
\end{proof}

Now, we can state the main result:
\begin{theorem}\label{theorem}
Let $\KK$ be a  characteristic zero field, and $V$ the vector space defined by $V=\KK x\oplus \KK y$. Let $F_0(x,y)$  be the Lie  series defined in definition \ref{def:FetG}.
On $\Lie(V)$, the Lie  series $F_0(x,y)$ and $G_0(x,y)=F_0(-y,-x)$ verify equation (\ref{KW1}) of the Kashiwara-Vergne conjecture:
$$x+y-\log(\exp(y)\exp(x))=(1-\exp(-\ad x))F(x,y)+(\exp (\ad y)-1)G(x,y)\ .$$
\end{theorem}
 
\begin{remark} Note that the theorem is true over any field $\KK$ of characteristic zero, i.e. such that $\KK$ contains $\QQ$, and not only for the two fields $\RR$ and $\CC$, as asked in the Kashiwara-Vergne conjecture, conjecture \ref{KV}.
\end{remark}

\begin{proof}
The proof of the theorem is done in three steps. First, we rewrite the Kashiwara-Vergne conjecture, thanks to the Baker-Campbell-Hausdorff series and its relation with the Eulerian idempotent (cf. proposition \ref{rewrite}). Secondly, we split the equation into two parts:
\begin{eqnarray*}
\Phi^-(x,y)=\Exp(-x)F(x,y)\ ,\\
\Phi^+(y,x)=-\Exp ( y)G(x,y)\ .
\end{eqnarray*}
By theorem \ref{splitKV} the first equation of this split equation admits $F_0(x,y)$ as a solution, and proposition \ref{prop:equationsplitsymmetry} gives $F_0(-y,-x)$ as a  solution of the second equation. 
Then adding the two equations and applying proposition \ref{cut}, proves the theorem.
\end{proof}

\begin{remark}\label{rem:Fexplicited}
The first terms of the symmetric solution $F_0(x,y)$ (cf. definition \ref{def:FetG}) are:

\begin{eqnarray*}
F_0(x,y)&=&\frac{1}{4}y+\frac{1}{24}xy-\frac{1}{24}yx\\
&&-\frac{1}{48}xxy+\frac{1}{24}xyx+\frac{1}{48}xyy-\frac{1}{48}yxx-\frac{1}{24}yxy+\frac{1}{48}yyx\\
&&-\frac{1}{180}xxxy+\frac{1}{60}xxyx+\frac{1}{480}xxyy-\frac{1}{60}xyxx-\frac{1}{240}xyxy\\
&&+\frac{1}{360}xyyy+\frac{1}{180}yxxx+\frac{1}{240}yxyx-\frac{1}{120}yxyy-\frac{1}{480}yyxx\\
&&+\frac{1}{120}yyxy-\frac{1}{360}yyyx\\
&&+\frac{1}{2880}xxxxy-\frac{1}{720}xxxyx-\frac{7}{2880}xxxyy+\frac{1}{480}xxyxx+\frac{7}{1440}xxyxy\\
&&+\frac{7}{2880}xxyyx+\frac{1}{720}xxyyy-\frac{1}{720}xyxxx-\frac{7}{720}xyxyx-\frac{1}{240}xyxyy\\
&&+\frac{7}{2880}xyyxx+\frac{1}{240}xyyxy-\frac{1}{360}xyyyx+\frac{1}{2880}yxxxx+\frac{7}{1440}yxyxx\\
&&+\frac{1}{240}yxyyx-\frac{1}{240}yyxyx-\frac{7}{2880}yyxxx+\frac{1}{720}yyyxx\\
&&+\textrm{higher order terms}
\end{eqnarray*}

This can be written (non-uniquely) in terms of bracket:
\begin{eqnarray*}
F_0(x,y)&=&\frac{1}{4}y+\frac{1}{24}[x,y]\\
&&-\frac{1}{48}(-[x,[x,y]]+[y,[y,x]])\\
&&-\frac{1}{180}[x,[x,[x,y]]]-\frac{1}{480}[x,[y,[x,y]]]-\frac{1}{360}[y,[y,[y,x]]]\\
&&-\frac{1}{240}[x,[x,[y,[x,y]]]]-\frac{7}{720}[x,[y,[x,[x,y]]]]+\frac{1}{144}[x,[y,[x,[x,y]]]]\\
&&-\frac{1}{240}[y,[x,[y,[x,y]]]]+\frac{1}{240}[y,[y,[x,[x,y]]]]+\textrm{higher order terms}
\end{eqnarray*}

\end{remark}

\section{Unicity of the solution}\label{sec:Unicity}

Moreover it can be proven that up to $\lambda x$, where $\lambda\in\KK$, the solution of the split equation (\ref{splitequation}) is unique.

\begin{proposition}
Let $F_0(x,y)$ be the Lie  series defined in definition \ref{def:FetG}.
 Any Lie  series $H(x,y)$ which is a solution of  (\ref{splitequation}):
 $$\Phi^-(y,x)=\Exp(- x)H(x,y) \ . $$
 is of the form $H(x,y)=F_0(x,y)+\lambda x$.
\end{proposition}

\begin{proof}
Let $H(x,y)$ be a Lie  series solution of the split equation (\ref{splitequation}). By proposition \ref{splitKV}, the  Lie series defined in definition \ref{def:FetG} is a solution of (\ref{splitequation}). Then substracting the two equations, it comes out that $H(x,y)-F_0(x,y)$ is in the kernel of $\Exp( (-x))$. Applying proposition \ref{Kernel} ends the proof.
\end{proof}

So we can conclude a unicity property for solutions of the split equation defined in the proof of theorem \ref{theorem}:

\begin{proposition}
Let $F_0(x,y)$ be the  Lie series defined in definition \ref{def:FetG}.
Let $(F(x,y),G(x,y))$ be a solution of equations (\ref{splitequation}):
\begin{eqnarray*}
\Phi^+(y,x)=-\Exp ( y)G(x,y)\ ,\\
\Phi^-(y,x)=\Exp(- x)F(x,y) \ . 
\end{eqnarray*}
Then, there exists $\lambda_1,\lambda_2\in\KK$ such that:
$F(x,y)=F_0(x,y)+\lambda_1$ and $G(x,y)=F_0(-y,-x)+\lambda_2$.

Conversely the pair $(F(x,y),G(x,y))$ is solution of the equations (\ref{splitequation}).\hfill$\Box$
\end{proposition}

Though the split equation (\ref{splitequation}) is unique up to the first degree term, the equation (\ref{KW1}) is not unique even if restricted to symmetric solutions, that is to say a pair solution $(F(x,y),F(-y,-x))$.
And we show that any non-symmetric solution of (\ref{KW1}) can be symmetrised thanks to the symmetry property verified by the left part of the Kashiwara-Vergne conjecture which is highlighted by Baker-Campbell-Hausdorff series. 

\begin{proposition}\label{BCHsym}The Baker-Campbell-Hausdorff series satisfy the following symmetry property:
$-\sum_{n\geq 2}\Phi_n(-y,-x)=\sum_{n\geq 2}\Phi_n(x,y)$.
\end{proposition}
\begin{proof}
By proposition \ref{cut}, the following identity holds:
$\sum_{n\geq 2}\Phi_n(x,y)=\Phi^+(x,y)+\Phi^-(x,y)$.
Moreover lemma \ref{phisymmetry} ensures that:
$\Phi^+(x,y)=-\Phi^-(-y,-x)\ .$
Therefore,
\begin{eqnarray*}
\sum_{n\geq 2}\Phi_n(x,y)&=&\Phi^+(x,y)+\Phi^-(x,y)\\
&=&-\Phi^-(-y,-x)-\Phi^+(-y,-x)\\
&=&-\sum_{n\geq 2}\Phi_n(-y,-x)\ .
\end{eqnarray*}
The proof is completed.
\end{proof}

From this proposition we deduce that any non-symmetric solution produces a symmetrised solution:

\begin{proposition}\label{symetrisation}
Let $(F(x,y),G(x,y))\in\Lie(V)^2$ be a non-symmetric solution of the first equation of the Kashiwara-Vergne conjecture (\ref{KW1}). The solution can be symmetrised in another solution:
\begin{eqnarray*}
F_1(x,y):=\frac{1}{2}(F(x,y)+G(-y,-x))+\lambda x \ ,\\
G_1(x,y):=\frac{1}{2}(G(x,y)+F(-y,-x))-\lambda y \ ,
\end{eqnarray*}
where $\lambda\in\KK$.
\end{proposition}

\begin{proof}
As $(F(x,y),G(x,y))$ is a solution of the first equation of the Kashiwara-Vergne conjecture (\ref{KW1}) it satisfies by proposition \ref{BCHsym}:
$$\sum_{n\geq 2}\Phi_n(y,x)=\Exp( (-x))F(x,y)-\Exp( y) G(x,y)\ ,$$
by exchanging $x$ and $-y$ we get:
$$\sum_{n\geq 2}-\Phi_n(-x,-y)=-\Exp( y)F(-y,-x)+\Exp( -x) G(-y,-x)\ .$$
We obtain the next equation by adding the two preceeding ones.
$$2\sum_{n\geq 2}\Phi_n(y,x)=\Exp( -x)(F(x,y)+G(-y,-x))-\Exp( y) (G(x,y)+F(-y,-x)) \ .$$
Moreover proposition \ref{Kernel}, permits the fact to add $\lambda x$ to the symmetrised solution $F(x,y)+G(-y,-x)$.
It is clear that the solution verifies the symmetry: $F_1(-y,-x)=G_1(x,y)$, which completes the proof.
\end{proof}

\section{Solution of the homogeneous equation}\label{section:sanssecondmembre}

We are interested in finding all solutions of the equation (\ref{KW1}) in the free Lie algebra generated by two indeterminates $x$ and $y$.
And therefore this section is devoted to solving the homogeneous equation  in order to set all the solutions of the equation (\ref{KW1}). That is to say solving the equation:
\begin{equation}\label{homogene}
\Exp(-x)F(x,y)=\Exp(y)(G(x,y)) \ ,
\end{equation}
for Lie polynomials $F(x,y)$ and $G(x,y)$.

\begin{lemma}\label{lemma:equationad}
The pair of Lie series $(P(x,y),Q(x,y))$ is a solution of the following equation:
\begin{equation}\label{eqavecadj}
\ad (x) P(x,y)+\ad (y) Q(x,y)=0\ ,
\end{equation}
if and only if there exists a non-commutative series $p(x,y)\in\Ker\gamma$ such that:
\begin{eqnarray*}
P(x,y)=\gamma((p(x,y))_x)\\
Q(x,y)=\gamma((p(x,y))_y)\ ,
\end{eqnarray*}
where $(p(x,y))_x$ (resp. $(p(x,y))_y$) denotes the x-part (resp. $y$-part) of $p(x,y)$.
\end{lemma}

\begin{proof}
Remark that the adjunction $ad(z)$, for $z\in V$, is a Lie homomorphism of degree $1$ and the Dynkin idempotent $\gamma$ is a degree preserving map. We can restrict the above equation to homogeneous polynomials $P(x,y)$ and $Q(x,y)$ of degree $n$.

Let $P(x,y)$ and $Q(x,y)$ be homogeneous Lie polynomials of degree $n$, defined as follows: 
\begin{eqnarray*}
P(x,y):=\gamma((p(x,y))_x)\\
Q(x,y):=\gamma((p(x,y))_y)\ .
\end{eqnarray*}
for a certain $p(x,y)\in\Ker\gamma_{n+1}$. Firstly, we verify that they satisfy equation (\ref{eqavecadj}).
Replacing $P$ and $Q$ in  equation (\ref{eqavecadj}) gives:
\begin{eqnarray*}
&&\ad(x)\gamma((p(x,y))_x)+\ad(y)\gamma((p(x,y))_y)\\
=&&\frac{n+1}{n}\gamma(x(p(x,y))_x+y(p(x,y))_y)\qquad\textrm{ (proposition \ref{gamma})}\\
=&&\frac{n+1}{n}\gamma(p(x,y))\qquad\textrm{ (definition \ref{def:x_part})}\\
=&&0\ ,
\end{eqnarray*}
as $p(x,y)\in\Ker\gamma_{n+1}$. Therefore the two defined Lie polynomials verify equation (\ref{eqavecadj}).

Conversely, let $P(x,y),Q(x,y)$ be homogeneous Lie polynomials of degree $n$. Suppose that the pair $(P,Q)$ verifies the equation:
\begin{equation*}
\ad (x) P(x,y)+\ad (y) Q(x,y)=0\ .
\end{equation*}
As $P(x,y)$ and $Q(x,y)$ are  Lie polynomials they verify the following property $\gamma(P(x,y))=P(x,y)$ and $\gamma(Q(x,y))=Q(x,y)$ respectively. 
So verifying the above equation is equivalent to verify:
\begin{equation}\label{equationavecgamma}
\ad (x) \gamma(P(x,y))+\ad (y) \gamma(Q(x,y))=0\ .
\end{equation} 
The linearity of $\gamma$ and proposition \ref{gamma} lead to:
\begin{equation*}
\gamma(xP(x,y)+yQ(x,y))=0\ .
\end{equation*}
And so $xP(x,y)+yQ(x,y)$ is a polynomial of degree $n+1$ which lies in the kernel of the Dynkin idempotent. We denote $p(x,y)$ this element of $\Ker\gamma_{n+1}$, where its $x$-part is $P(x,y)=(p(x,y))_x$ and its $y$-part is $Q(x,y)=(p(x,y))_y$. 
Recall that  $P(x,y)$ and $Q(x,y)$ are  Lie series verifying $P(x,y)=\gamma(P(x,y))=\gamma((p(x,y))_x)$ (resp. $Q(x,y)=\gamma(Q(x,y))=\gamma((p(x,y))_y)$).
Therefore $P(x,y)=\gamma((p(x,y))_x)$ and $Q(x,y)=\gamma((p(x,y))_y)$.

The proof is completed.
\end{proof}

\begin{example}\label{ex:Vergne}
M. Vergne found that the polynomial $$P(x,y)=[x,[y,[x,[x,y]]]]-2[y,[x,[x,[x,y]]]-[y,[y,[y,[y,x]]\neq 0$$
 verifies the equation (\ref{eqavecadj}): $$ad(x)P(x,y)+ad(y)P(-y,-x)=0\ ,$$ (private communication). Under lemma \ref{lemma:equationad} we should be able to prove that the polynomial $p(x,y):=xP(x,y)+yP(-y,-x)$ is in the kernel of the Dynkin idempotent.  
This is true as there exists a polynomial:
\begin{eqnarray*}
q(x,y)&=&2xxxxyy-8xxxyxy+xxxyyx+12xxyxxy\\
&&-4xxyxyx+xxyyxx-2xxyyyy-8xyxxxxy\\
&&+6xyxxyx-4xyxyxx+xyyxxx+8xyxyyy\\
&&-12xyyxyy+8xyyyxy-xyyyyx\\
&&+yxxxxy-yxxyyy+4yxyxyy-6yxyyxy-yyxxyy\\
&&+4yyxyxy-yyyxxy
\end{eqnarray*} 
is such that $p(x,y)=q(x,y)-\gamma(q(x,y))$ which is in the kernel of the Dynkin idempotent.
\end{example}

\begin{proposition}\label{prop:lessolutionshomogenes}
 The pair  $(F(x,y),G(x,y))$ is  a solution of the homogeneous equation (\ref{homogene}) if and only if there exists an element $p(x,y)$ in the kernel the Dynkin idempotent such that  $F(x,y)$ is equal to $$F(x,y):=\B (-x)\gamma((p(x,y))_x)+\lambda_1 x$$ and $G(x,y)$ is equal to $$G(x,y):=\B (y)\gamma((p(x,y))_y)+\lambda_2 y\ ,$$ where $(p(x,y))_x$ (resp. $(p(x,y))_y$) denotes the x-part (resp. $y$-part) of $p(x,y)$ and $\lambda_1,\lambda_2\in\KK$.
\end{proposition}

\begin{proof}
First we verify that the Lie series $F(x,y)$ and $G(x,y)$ are solutions of equation (\ref{homogene}).
As the kernel of the Dynkin idempotent can be described homogeneously we can restrict ourselves to elements of the kernel of $\gamma_n:\Vt n\rightarrow \Vt n$. 

Let $p_n(x,y)\in\Ker\gamma_n$, define the two Lie polynomials $F(x,y)$ and $G(x,y)$ as follows $F(x,y):=\B(-x)\gamma((p_n(x,y))_x)$ respectively $G(x,y):=\B((y))\gamma((p_n(x,y))_y)$.
Then the equation  
$\Exp(-x)F(x,y)-\Exp((y))(G(x,y))$
becomes by definition of $F(x,y)$ and $G(x,y)$ and proposition \ref{prop:inverse}:
$$ad(-x)\gamma((p_n(x,y))_x)-\ad(y)\gamma((p_n(x,y))_y)$$
Lemma \ref{lemma:equationad} assures that the pair $(F(x,y),G(x,y))$ is solution of the homogeneous equation.

Let $F(x,y)$ and $G(x,y)$ be a Lie series which are  solutions of equation (\ref{homogene}).
By proposition \ref{Image}, there exists two Lie series $P(x,y)$ and $Q(x,y)$ such that  $\Exp(  -x)F(x,y)=\ad (-x)P(x,y)$ and $\Exp(  y)G(x,y)=\ad (-x)Q(x,y)$ respectively. 
As $(F(x,y),G(x,y))$ is a solution of  equation (\ref{homogene}) the pair $(P(x,y),Q(x,y))$ is a solution of the following equation:
\begin{equation*}
\ad (x) P(x,y)+\ad (y) Q(x,y)=0\ .
\end{equation*}
And we note $$p(x,y):=xP(x,y)+yQ(x,y)\ .$$
Therefore lemma \ref{lemma:equationad} gives that $P(x,y)=\gamma((p(x,y))_x)$ and $Q(x,y)=\gamma((p(x,y))_y)$.
As the pair $(F(x,y),P(x,y))$ verify the equation (\ref{eq:passage}) we can 
apply lemma \ref{lemma:passeadexp} to obtain:
$$F(x,y)=\B( (-x))P(x,y)+\lambda_1 x\ ,$$
for a certain $\lambda_1\in\KK$. An analogue to the lemma \ref{lemma:passeadexp} (changing $\Exp(-x)$ in $\Exp(y)$) gives that:
$$G(x,y)=\B(y)Q(x,y)+\lambda_2 y\ ,$$
for a certain $\lambda_1\in\KK$. 
This ends the proof.
\end{proof}

\medskip
We would like to have a more explicit formula for the Lie series solution of the homogeneous equation (\ref{homogene}). So we construct maps $\Psi_{x}:T(V)\rightarrow T(V)$ which will simplify the span of the vector space of solutions.
\begin{definition}\label{Psi}
We denote by $\Psi_{x}:T(V)\rightarrow T(V)$ the map defined by $p(x,y)\mapsto\gamma((p(x,y)-\gamma(p(x,y))_x))$.
\end{definition}

This map can be explicited in terms of permutations as follows:
\begin{proposition}
The map $\Psi_{x}:T(V)\rightarrow T(V)$
is induced by 
\begin{eqnarray*}
&&x_{1}\ldots x_{n}\mapsto \frac{1}{n-1}\big((-1)^{n-2}\delta_{x,x_j}
(\sum_{\begin{array}{c}
\scriptstyle 1\leq k\leq n-2\\
\scriptstyle \tau\in D_{\{1,\cdots k\}\subset S_{n-1}}
\end{array}}
(n-1)x_{\tau(2)}\ldots x_{\tau(n)})\\
&&+\sum_{\begin{array}{c}
\scriptstyle 1\leq k\leq n-2\\
\scriptstyle 1\leq j\leq n-2
\end{array}} 
\sum_{\begin{array}{c}
\scriptstyle \sigma\in D_{\{1,\cdots k\}}\subset S_n\\
\scriptstyle \omega\in D_{\{1,\cdots j\}}\subset S_{n-1}
\end{array}}
(-1)^{k+j-1}
\delta_{x,x_{\sigma(n)}}(x_{\omega(\sigma(n-1))}\ldots x_{\omega(\sigma(1))})\big)\ ,
\end{eqnarray*}
where the map $\delta_{x,x_j}:T(V)\rightarrow \KK$ is the map induced by:
\begin{eqnarray*}
x_{1}\ldots x_{n}\mapsto
\left\{
\begin{array}{rl}
1 & \mbox{if } x_{j} = x \\
0 & \mbox{otherwise }
\end{array}
\right.
\end{eqnarray*}
\end{proposition}
\begin{proof}
Define the non-commutative polynomial $p_n(x_1\ldots x_n):=n(x_1\ldots x_n-\gamma(x_1\ldots x_n))$
which can be explicited by proposition \ref{prop:KernelDynkin} as: $$p_n(x_1\ldots x_n):=(n-1)x_1\ldots x_n+\sum_{k=1}^{n-2}\sum_{\sigma\in D_{\{1,\cdots k\}}}(-1)^{n+k-1} (x_n\ldots x_1)^\sigma\ ,$$ for $x_i\in V$ and $n\geq 2$.

Let $x_1\ldots x_n$ be a monomial of degree $n$ in $\Vt n$. Then, the $x$-part of this monomial is the monomial $x_2\ldots x_n$ if $x_1=x$. This can be sum up as $(x_1\ldots x_n)_x=\delta_{x,x_j}x_2\ldots x_n$.

Therefore the $x$-part of $p_n(x_1\ldots x_n)$ is 
$$
(p_n(x_1\ldots x_n))_x=(n-1)\delta_{x,x_j}x_2\ldots x_n+\sum_{k=1}^{n-2}\sum_{\sigma\in D_{\{1,\cdots k\}}}(-1)^{n+k-1} \delta_{x,x_{\sigma(n)}}(x_{\sigma(n-1)}\ldots x_{\sigma(1)})\ .
$$ 
By the explicit formula of proposition \ref{explicitDynkin}:
$$
\gamma_n(x_1\ldots x_n)=\frac{(-1)^{n-1}}{n}\sum_{k=0}^{n-1}(-1)^k\sum_{\sigma\in D_{\{1,\ldots,k\}}}(x_n\ldots x_1)^\sigma \ ,
$$ 
taking the Dynkin idempotent of $(p_n(x_1\ldots x_n))_x$
gives the following:
\begin{eqnarray*}
&&\gamma((p_n(x_1\ldots x_n))_x)=(n-1)\delta_{x,x_j}\gamma(x_2\ldots x_n)\\
&&\qquad+\sum_{\begin{array}{c}
\scriptstyle 1\leq k\leq n-2\\
\scriptstyle \sigma\in D_{\{1,\cdots k\}\subset S_{n}}
\end{array}}
(-1)^{n+k-1} \delta_{x,x_{\sigma(n)}}\gamma(x_{\sigma(n-1)}\ldots x_{\sigma(1)})\\
&&\qquad=(-1)^{n-2}\delta_{x,x_j}
(\sum_{\begin{array}{c}
\scriptstyle 1\leq k\leq n-2\\
\scriptstyle \tau\in D_{\{1,\cdots k\}}\subset S_{n-1}
\end{array}}
x_{\tau(2)}\ldots x_{\tau(n)})\\
&&\qquad+(-1)^{n-2}\sum_{\begin{array}{c}
\scriptstyle 1\leq k\leq n-2\\
\scriptstyle \sigma\in D_{\{1,\cdots k\}}
\end{array}} \delta_{x,x_{\sigma(n)}}(
\sum_{\begin{array}{c}
\scriptstyle 1\leq j\leq n-2\\
\scriptstyle \omega\in D_{\{1,\cdots j\}\subset S_{n-1}}
\end{array}}
\frac{(-1)^{n+k+j-1}}{n-1}
x_{\omega(\sigma(n-1))}\ldots x_{\omega(\sigma(1))})\\
&&=n\Psi_{x}(x_1\ldots x_n)
\end{eqnarray*} 
\end{proof}

\begin{proposition}
Let $V$ be the vector space spanned by the indeterminates $x$ and $y$.
Let $\Psi_{x},\Psi_{y}$ be the maps defined in the above definition (definition \ref{Psi}).
Let $(F(x,y),G(x,y))\in\Lie(V)^2$ be a solution of the homogeneous equation (\ref{homogene}).
Then there exists a polynomial $m(x,y)\in T(V)$ such that:
\begin{eqnarray*}
F(x,y)=\B( -x)\Psi_{x}(m(x,y))\ ,\\
G(x,y)=\B( y)\Psi_{y}(m(x,y))\ .
\end{eqnarray*}
\end{proposition}
\begin{proof}
Recall that the kernel of the Dynkin idempotent is generated by $1$ and elements $$p_n(x_1\ldots x_n):=(n-1)x_1\ldots x_n+\sum_{k=1}^{n-2}\sum_{\sigma\in D_{\{1,\cdots k\}}}(-1)^{n+k-1} (x_n\ldots x_1)^\sigma\ ,$$ for $x_i\in V$ and $n\geq 2$.Then, proposition \ref{prop:lessolutionshomogenes} and proposition \ref{prop:KernelDynkin} complete the proof.
\end{proof}

\begin{example}With the notations of  example \ref{ex:Vergne}
we know that $$(\B(-x)P(x,y),\B(y)P(-y,-x))$$ is a  solution of the homogeneous equation (\ref{homogene}). This solution can be explicited as in the above proposition since: 
$$P(x,y)=\gamma(P(x,y))=\gamma((q(x,y)-\gamma(q(x,y)))_x)=\Psi_x(q(x,y))\ .$$
\end{example}

We can now state the theorem giving all solutions of equation (\ref{KW1}) in the free Lie algebra generated by the two non-commutative indeterminates $x$ and $y$:

\begin{theorem}\label{theorem:allsolution}
Let $V$ be the $\KK$-vector space spanned by the indeterminates $x$ and $y$.
Let $(F_0(x,y),F_0(-y,-x))$ be the particular symmetric solution of equation (\ref{KW1}) constructed in definition \ref{def:FetG}.
Let $\Psi_{x},\Psi_{y}$ be the maps defined in definition \ref{Psi}.

Let the pair $(F(x,y),G(x,y))\in\Lie(V)^2$ of Lie  series be a solution of the equation (\ref{KW1}). Then there exists polynomial $p(x,y)\in T(V)$ such that:
\begin{eqnarray*}
F(x,y)=F_0(x,y)+\sum_{m\geq 1}\frac{B_m}{m!}(-1)^m(\ad(x))^{ m}
\Psi_{x}(p(x,y))\\
G(x,y)=F_0(-y,-x)+\sum_{m\geq 1}\frac{B_m}{m!}(\ad(y))^{ m}\Psi_{y}(p(x,y))
\end{eqnarray*}
Conversely, the pair $(F(x,y),G(x,y))$ is a solution of equation (\ref{KW1}).
\end{theorem}

\begin{proof}
It is clear by proposition \ref{prop:lessolutionshomogenes} and by theorem \ref{theorem} that the pair $(F(x,y),G(x,y))$ is solution of the equation (\ref{KW1}).
Conversely, let $(F(x,y),G(x,y))$ be a solution of equation (\ref{KW1}). As $(F_0(x,y),F_0(-y,-x))$ is also a solution, their difference is solution of the homogeneous equation (\ref{homogene}). And proposition \ref{prop:lessolutionshomogenes} ends the proof. 
\end{proof}

\section{Another decription of all the solutions of equation (\ref{KW1})}\label{sec:otherview}

In \cite{PR} and in \cite{C} there is another description of the kernel of the Dynkin idempotent recalled as proposition \ref{kernelgamma}. 
The same proof as before brings us to another formulation of theorem \ref{theorem:allsolution}. In order to simplify the statement of this theorem we introduce a few notations.

\begin{definition}\label{def:EetA}
Let $p(x,y)$ be a non-commutative polynomial in indeterminates $x$ and $y$.
We define:
\begin{eqnarray*}
A(p)(x,y)&:=&\gamma(p(x,y))p(x,y)-\gamma(p(-y,-x))p(-y,-x)\\
E(p)(x,y)&:=&\gamma(p(x,y))p(x,y) \ ,
\end{eqnarray*} 
Moreover we define $A(p)_n$ (resp. $E(p)_n$) as the homogeneous part of $A(p)$ (resp. $E(p)$) of degree $n$, as a Lie polynomial can uniquely be seen as a non-commutative polynomial. 
\end{definition}

\begin{proposition}
Let $V$ be the vector space spanned by the indeterminates $x$ and $y$.
Let $(F_0(x,y),F_0(-y,-x))$ be the particular symmetric solution of equation (\ref{KW1}) constructed in definition \ref{def:FetG}.
Let the pair  $(F(x,y),G(x,y))\in\Lie(V)^2$ of Lie  series in indeterminates $x$ and $y$.

If $(F(x,y),G(x,y))$ is a solution of the equation (\ref{KW1}) then, there exists $\lambda_1,\lambda_2\in\KK$  and a finite set of indices $I$ indexing a finite family of non-commutative polynomials $p^i(x,y)\in T(V)$  and a finite family of scalars $\mu_i\in\KK$ such that:
\begin{eqnarray*}
F(x,y):= F_0(x,y) +\sum_{i\in I}\B(-x)\gamma(\sum_{n\geq 0}\frac{n}{n+1}\mu_i (E(p^i)_n(x,y))_x)+\lambda_1 x \ ,\\
G(x,y):= F_0(-y,-x) +\sum_{i \in I}\B(y)\gamma(\sum_{n\geq 0}\frac{n}{n+1}\mu_i (E(p^{i})_{n}(x,y))_y)+\lambda_2 y \ ,
\end{eqnarray*}
with the notations of definition \ref{def:EetA}.

Conversely, the pair $(F(x,y),G(x,y))$ is a solution of equation (\ref{KW1}).
\end{proposition}

With this point of view there is a way of producing all symmetric solutions of equation (\ref{KW1}) as follows:

\begin{proposition}
Let $V$ be the vector space spanned by the indeterminates $x$ and $y$.
Let $(F_0(x,y),F_0(-y,-x))$ be the particular symmetric solution of equation (\ref{KW1}) constructed in definition \ref{def:FetG}.
Let $F(x,y)\in \Lie(V)$ be a Lie  series in indeterminates $x$ and $y$.

If $(F(x,y),F(-y,-x))$ is a symmetric solution of the equation (\ref{KW1}) then, there exists $\lambda_1\in\KK$  and a finite set of indices $I$ indexing a finite family of non-commutative polynomials $p^i(x,y)\in T(V)$  and a finite family of scalars $\mu_i\in\KK$ such that:
\begin{eqnarray*}
F(x,y):= F_0(x,y) +\sum_{i \in I}\B(-x)\gamma(\sum_{n\geq 0}\frac{n}{n+1}\mu_i (A(p^i)_n(x,y))_x)+\lambda_1 x \ ,
\end{eqnarray*}
with the notations of definition \ref{def:EetA}.

Conversely, the pair $(F(x,y),F(-y,-x))$ is a symmetric solution of equation (\ref{KW1}).
\end{proposition}

The proof is analogous to the one used for proposition \ref{prop:lessolutionshomogenes}. It is to be noted that instead of elements of the kernel of the Dynkin idempotent we need the anti-symmetric elements of the kernel of the Dynkin idempotent which is spanned by:

\begin{proposition}\label{prop:KernelDynkinantisymmetric}
The anti-symmetric elements of the kernel of the Dynkin idempotent are spanned by the elements $$A(p)(x,y)=\gamma(p(x,y))p(x,y)-\gamma(p(-y,-x))p(-y,-x)\ ,$$ for $p(x,y)\in T(V)$.
\end{proposition}

\begin{proof}
Let $p(x,y)$ denote a non-commutative polynomial.
By proposition  \ref{kernelgamma} it is clear that the elements spanned by $A(p)(x,y)$ are anti-symmetric elements of the kernel of the Dynkin idempotent.

Conversely, let $q(x,y)$ be an anti-symmetric element of the kernel of the Dynkin idempotent. By proposition \ref{kernelgamma} there exist a finite family $p^i(x,y)\in T(V)$ and $\lambda_i\in\KK$ such that:  $$q(x,y)=\sum_{i\geq 0}\lambda_i\gamma(p^i(x,y))p^i(x,y)\ .$$
 By the anti-symmetry property of $q(x,y)$ this sum can also be rewritten as $$q(x,y)=\frac{1}{2}\sum_{i\geq 0}\lambda_i(\gamma(p^i(x,y))p^i(x,y)-\gamma(p^i(-y,-x))p^i(-y,-x))\ .$$
Therefore anti-symmetric elements of the kernel of the Dynkin idempotent  are spanned by elements $\gamma(p(x,y))p(x,y)-\gamma(p(-y,-x))p(-y,-x)$, where $p(x,y)\in T(V)$, which ends the proof.
\end{proof}

\section{Multilinearised Kashiwara-Vergne conjecture}\label{sec:multilinear}

Remark that the previous method can be extended in order to find all the solutions in $\Lie(\KK x_1\oplus\cdots\oplus\KK x_n)$ of the first equation of the following multilinear version of Kashiwara-Vergne conjecture:
\begin{conjecture}
For any Lie algebra $\gg$ of finite dimension, we can find  series $F_1,\ldots, F_n$ such that they satisfy
\begin{enumerate}
\addtocounter{enumi}{\value{equation}}
\item \begin{eqnarray*}
x_1+\cdots+x_n-\log(e^{x_n}\cdots e^{x_1})=(1-e^{-\ad x_1})F_1(x_1,\ldots,x_n)\\
+(1-e^{\ad x_2})F_2(x_1,\ldots,x_n)+\cdots+(1-e^{(-1)^n \ad  x_n})F_n(x_1,\ldots,x_n)\ ,
\end{eqnarray*}\label{KW1multilinear}
\item $F_1,\ldots,F_n$ give $\gg$-valued convergent power series on $(x_1,\ldots,x_n)\in \gg^{\times n}\ ,$\label{KW2multilinear}
\item \label{KW3multilinear}\begin{eqnarray*}
&&\tr(\ad x_1\circ\partial_{x_1} F_1;\gg)+\cdots+\tr(\ad x_n\circ\partial_{x_n} F_n;\gg)=\\
&&\frac{1}{2}\tr(\frac{\ad x_1}{e^{\ad x_1} -1}+\cdots+\frac{\ad x_n}{e^{\ad x_n} -1}+\frac{\ad \Phi(x_n,\ldots,x_1)}{e^{\ad \Phi(x_n,\ldots,x_1)} -1}-1;\gg)\ .
\end{eqnarray*}
\end{enumerate}
Here $\Phi(x_1,\ldots,x_n)=\log(e^{x_1}\cdots e^{x_n})$ and $\partial_{x_i} Fi$  is the $\mathbf{End}(\gg)$-valued real analytic function defined by
\begin{equation*}
\gg\ni a\mapsto \frac{d}{dt}F(x_1,\ldots,x_{i-1},x_i+ta,x_{i+1},\ldots,x_n)|_{t=0} 
\end{equation*}
and $\tr$ denotes the trace of an endomorphism of $\gg$.
\end{conjecture} 
We are able to construct explicitly a particular solution $F_{1,0},\ldots,F_{n,0}$ of  equation (\ref{KW1multilinear}).

\begin{theorem}\label{theo:particularmultilinear}
Let $\KK$ be a  characteristic zero field, and $V$ the vector space defined by $V=\KK x_1\oplus\cdots\oplus \KK x_n$. Let $F_{i,0}(x,y)$ be the Lie  series defined below:
\begin{eqnarray*} 
a_i(x_1,\ldots,x_n)&:=&\sum_{m\geq 1}\frac{m}{m+1}\sum_{\begin{array}{c}
\scriptstyle i_1+\cdots+i_n=m \\
\scriptstyle i_1,\ldots,i_n\geq 1
\end{array}} \frac{1}{i_1!}\cdots\frac{1}{(i_k+1)!}\cdots\frac{1}{i_n!}\\
&&\gamma\circ (e_n(\underbrace{x_n,\ldots,x_n}_{i_n},\ldots,\underbrace{x_1,\ldots,x_1}_{i_1}))_{x_i} \ .
\end{eqnarray*}
\begin{eqnarray*}
&&F_{i,0}(x_1,\ldots,x_n):=-\sum_{m\geq 0}\frac{B_m}{m!}(-1)^m(\ad x_i)^{m}\circ a_i(x_1,\ldots,x_n)\ .
\end{eqnarray*}
On $\Lie(V)$, the Lie  series $F_{1,0}(x_1,\ldots,x_n),\ldots,F_{n,0}(x_1,\ldots,x_n)$ verify equation (\ref{KW1multilinear}) of the multilinear Kashiwara-Vergne conjecture:
\begin{eqnarray*}
x_1+\cdots+x_n-\log(e^{x_n}\cdots e^{x_1})=(1-e^{-\ad x_1})F_1(x_1,\ldots,x_n)\\
+(1-e^{\ad x_2})F_2(x_1,\ldots,x_n)+\cdots+(1-e^{(-1)^n \ad  x_n})F_n(x_1,\ldots,x_n)
\end{eqnarray*}
\end{theorem}

Remark that we can also give an analogous to proposition \ref{splitKV} where the solutions will be unique up to $(\lambda_1 x_1,\ldots,\lambda_n x_n)$ for $\lambda_i\in\KK$. 

We can state the theorem giving explicitly all the solutions of  equation (\ref{KW1multilinear}):

\begin{theorem}\label{theorem:allsolutionmultilinear}
Let $V$ be the vector space spanned by the indeterminates $x_1,\ldots,x_n$.
Let  $F_{1,0}(x_1,\ldots,x_n),\ldots,F_{n,0}(x_1,\ldots,x_n)$ be the particular symmetric solution of equation (\ref{KW1}) constructed in theorem \ref{theo:particularmultilinear}.
Let $\Psi_{x_i}:T(V)\rightarrow T(V)$ be the map defined in definition \ref{Psi}.

If the pair $(F_{1}(x_1,\ldots,x_n),\ldots,F_{n}(x_1,\ldots,x_n))\in\Lie(V)^n$ of Lie  series in indeterminates $x_1,\ldots,x_n$ is a solution of the equation (\ref{KW1multilinear}), then there exists a polynomial $m(x,y)$ such that:
\begin{eqnarray*}
F_i(x,y)=F_{i,0}(x_1,\ldots,x_n)+\B((-1)^ix_i)\Psi_{x_i}(m(x,y))\\
\end{eqnarray*}
Conversely, the n-tuple $(F_1,\ldots,F_n)$ is a solution of equation (\ref{KW1multilinear}).
\end{theorem}

The above theorem has an analogous version with  the  description of the kernel of the Dynkin idempotent due to Patras and Reutenauer's:

\begin{proposition}
Let $V$ be the vector space spanned by the indeterminates $x_1,\ldots,x_n$.
Let $F_{1,0}(x_1,\ldots,x_n),\ldots,F_{n,0}(x_1,\ldots,x_n)$ be the particular solution of equation (\ref{KW1multilinear}) constructed in theorem \ref{theo:particularmultilinear}.
Let the n-tuple  $(F_{1}(x_1,\ldots,x_n),\ldots,F_{n}(x_1,\ldots,x_n))\in\Lie(V)^n$ of Lie  series in indeterminates $x_1,\ldots,x_n$.

If $(F_{1},\ldots,F_{n})$ is a solution of the equation (\ref{KW1}) then, there exists there exists $\lambda_1,\ldots,\lambda_n\in\KK$  and a finite set of indices $J$ indexing a finite family of non-commutative polynomials $p^j(x_1,\ldots,x_n)\in T(V)$  and a finite family of scalars $\mu_j\in\KK$ such that:
\begin{eqnarray*}
&&F_{i}(x_1,\ldots,x_n):= F_{i,0}(x_1,\ldots,x_n) \\
&&\qquad\qquad+\B((-1)^ix_i)\gamma\big(\sum_{j\in J}\sum_{m\geq 0}\frac{m}{m+1}\mu_j (E(p^j)_n(x_1,\ldots,x_n))_{x_i}\big)+\lambda_{i} x_i \ ,\\
\end{eqnarray*}
with the notations of definition \ref{def:EetA}.

Conversely, the $n$-tuple is solution of equation (\ref{KW1multilinear}).
\end{proposition}
Indeed, there exists an analogous formula linking the multilinearized Baker-Campbell-Hausdorff formula defined as $\Phi(x_1,\ldots,x_n):=\log(e^{x_1}\cdots e^{x_n})$ and the Eulerian idempotent (cf.~\cite{L}).
All the proofs will be analogous as the Eulerian and the Dynkin idempotent are defined on the tensor module $T(V)$ over any vector space $V$ and can be particularized in the case where $V=\KK x_1\oplus\cdots\oplus\KK x_n$.

The case treated in this paper is the case where $n=2$, $F_1(x,y)=F(x,y)$ and $F_2(x,y)=-G(x,y)$.

\end{document}